\documentclass[a4paper,12pt]{article}
\usepackage{amsmath,graphicx}
\usepackage{amssymb}

\newcommand{\be}{\begin{equation}}
\newcommand{\ee}{\end{equation}}
\newcommand{\ba}{\begin{eqnarray}}
\newcommand{\ea}{\end{eqnarray}}
\newcommand{\bee}{\begin{equation*}}
\newcommand{\eee}{\end{equation*}}
\newcommand{\baa}{\begin{eqnarray*}}
\newcommand{\eaa}{\end{eqnarray*}}
\begin{document}

\begin{titlepage}

\begin{center}

\includegraphics[width=0.15\textwidth]{./ISI}\\[1cm]

\textsc{\LARGE Indian Statistical Institute}\\[1.5cm]

\textsc{\Large M.Stat II year Project Report}\\[0.5cm]

\hrule
\hspace*{\fill} \\[4mm]

{ \huge \bfseries A Recurrent Rotor-Router Configuration in $\mathbb{Z}^3$}\\[0.4cm]

\hrule
\hspace*{\fill} \\[4mm]

\begin{minipage}{0.4\textwidth}
\begin{flushleft} \large
\emph{Student:}\\
Tulasi Ram Reddy A \textsc{\emph{MB-0832}}
\end{flushleft}
\end{minipage}
\begin{minipage}{0.4\textwidth}
\begin{flushright} \large
\emph{Supervisor:} \\
Dr. Arni S.R. Srinivasa Rao \textsc{}
\end{flushright}
\end{minipage}

\vfill

{\large \today}

\end{center}

\end{titlepage}

\section{Introduction to Rotor-Router Model}
\paragraph{}
Informally a random walk is a process in which a particle takes a successive random steps. A pseudo random walk appears to be a random walk, but it is not. It exhibits randomness in some sense but is obtained from a deterministic procedure. One of the popular random walk is the one that is executed on the lattice points. In the further paragraphs you will be able to notice that, the walks executed by the particles in Rotor-Router Model together can be viewed as a pseudo random walk.
\paragraph{}
The need to develop the Rotor-Router Model arose during the study of Internal Diffusion-limited Aggregation (IDLA) and abelian sandpile models. IDLA is a model in which a particle executes a random walk in the lattice \textbf{$\mathbb{Z}^{d}$} beginning at the origin and ends at the point when it first reaches a point outside the set \textbf{$A$} which contains origin. Now the end point of this random walk is adjoined to the set \textbf{$A$}, and this process is iterated. The fact that after \emph{n} walks have been executed, the resulting region \textbf{$A$}
rescaled by a factor  $n^{1/d}$ can be approximated by an Euclidian sphere in $\mathbb{R}^d$ as $\emph{n} \longrightarrow \infty$ was proved by Lawler(\cite{LBG}, 1992). Whereas in an abelian sandpile model, introduced by Bak et al.(\cite{BTK}, 1988) and studied by Dhar(\cite{DD}, 1990) and Bonjour, et al.(\cite{BLP}, 1991), each of the points in \textbf{$A$} were assigned a non-negative integer, which represent the number of grains of sand at that point. Origin is assigned with a single grain of sand and at each step, every point that has at least $2d$ grains, ejects a grain to each of its neighbours. The process is repeated until every point contains at most $2d-1$ grains. Many of the fundamental mathematical properties which remain conjectural in IDLA can actually be solved in abelian sandpile model. Whereas no asymptotic results were known for sandpile models and it also appears likely that its asymptotics are not spherical. Though they have a common abelian property, there are considerable differences between IDLA and sandpile model. So a need arose to develop a model that is intermediate to both.
\paragraph{}
Rotor-Router Model introduced by Jim Propp(\cite{LLRR}) can be visualized as a deterministic analogue of IDLA. First a cyclic ordering of the $2$ directions that are along the co-ordinate axes in $\mathbb{Z}^d$ is fixed. A 'rotor' is placed at each point in \textbf{$A$} which specifies one of the $2d$ directions. When a particle is passed through a rotor it  changes its direction to the direction that is following immediately after its previous direction in the cyclic ordering. A particle starts at origin and is routed through these points until it reaches a point that is not in the initial set \textbf{$A$}. Once the particle reaches to a point outside the set $A$, the same point is embedded to this set placing a rotor having an initial direction same as that of the origin.Then a new particle will be starting from the origin and follows the same procedure in new larger set. The process will be continued in this manner.
\paragraph{}
These kind of models are very useful for computational purposes. They converge much faster than the conventional random processes.

\section{Rotor walk as a substitute to Random walk}
\paragraph{}
Consider the sequence $\{1,0,1,0,1,0 ,....\}$. This can be seen as a realization of the sequence of \emph{iid} random variables $\{X_i\}$ which takes values 0 or 1 with equal probability. Now consider the case of simple random walk in $\mathbb{Z}^d$. Suppose if we have a collection of sequences of \emph{iid} random variables at each lattice point, where each random variable assumes all the 2d directions with same probability. Now start at the origin and execute the random walk in the following manner. Whenever the particle reaches a point in $\mathbb{Z}^d$ for the $n^{th}$ time, look at the $n^{th}$ random variable of the corresponding sequence present at that point. Based on the direction obtained from that random variable make the next move. In this way a simple random walk on $\mathbb{Z}^d$ is executed.
\paragraph{}
Now instead of the sequences $\{X_i\}$, one can use the sequences of the form $\{e_0,e_1,..,e_{2d-1},e_0,e_1,..,e_{2d-1},e_0,e_2,..,e_{2d-1},.....\}$, where $\{e_0,e_1,....,e_{2d-1}\}$ are the natural 2d directions in $\mathbb{Z}^d$. These sequences can be viewed as a realizations for the sequences $\{X_i\}$. The walk that is obtained now is a rotor walk. In this way one can see the rotor walk on $\mathbb{Z}^d$ as a pseudo random walk on $\mathbb{Z}^d$. Rotor walk in general exhibits the properties of the random walk.

\section{Mathematical formulation of Rotor Walk}
\paragraph{}
Let $E_d = \{{e_0,...,e_{2d-1}}\}$, be the set of $2d$ directions that are along the co-ordinate axes in $\mathbb{Z}^d$ in that given order. A state of a Rotor-Router process is represented by a pair $(x,l)$, where $x \in \mathbb{Z}^d$ represents the location of the particle, and $l : \mathbb{Z}^d \longrightarrow \{0,1,...2d-1\}$ indicates the direction of the rotor at each point. Define $g(x,l) = (x+e_{l(x)},l_x)$, where $l(x)$ is the labeling given by

\be
l_x(x') = \left\{ \begin{array}{cc} \displaystyle l(x')$ \emph{+ 1(mod 2n)} $ $  \emph{if} $ x'=x; \\l(x'), $             \emph{ otherwise.}$ \end{array}\right.
\ee

This explains the position to which the particle is moved and also the direction to which the rotor at that point is shifted. Composing the function $g$ to itself at the state $(0,l)$ we can obtain the path of the particle as $(0,x_1,x_2,...)$ is a lattice path in $\mathbb{Z}^d$, beginning at the origin. We denote this path by $p=p(l)$. Earlier we have assumed that the particle reaches outside the finite region containing the initial point with in finitely many steps. It will be proved in the following lemma.
\paragraph{}
\textbf{Lemma:} \emph{Let $A$ be the finite set containing the origin. Then the lattice path $p(l)$ leaves the region $A$ in finitely many steps.} \\
\textbf{Proof:} Suppose not, then $\exists$ a lattice point in $A$ such that the particle is visited infinitely many times. And for any point which is visited infinitely many times, the same true for its neighbouring lattice points.  So any point that is connected with the point that is visited infinitely many times is also visited infinitely many times. Inducting the above argument along the path from origin to the point outside $A$, we obtain that the particle will reached outside the region $A$. But it contradicts our assumption. Hence proved.
\paragraph{}
Thus the entire process in the Rotor-Router Walk can be interpreted mathematically by the above terminology used by Levine(\cite{LLRR}).

\subsection*{Recurrent States:}
\paragraph{}
We say a state in $\mathbb{Z}^d$ is `\emph{recurrent}', if a particle started at origin visits that state infinitely often. It is easy to see that if a state in $\mathbb{Z}^d$ is recurrent then every other state in $\mathbb{Z}^d$ is also recurrent. So the initial rotor configuration for which every state is recurrent is called \emph{`Recurrent Rotor Configuration'}. So if a state is not recurrent then we call it a \emph{`Transient state'},i.e. a particle starting at origin visits this state only finitely many times, so is every other state. Hence the initial rotor configuration which is not recurrent is a \emph{`Transient Rotor Configuration'}\\
\subsubsection*{An open problem from Holroyd and Propp(\cite{HOPR}) }
\paragraph{}
We know that in a classical random walk, every state in $\mathbb{Z}^d$ is recurrent for $d\leq2$ and is transient for $d\geq3$. We can easily find a several examples of initial rotor configurations for which it is recurrent in $\mathbb{Z}^d(d\leq2)$, transient in $\mathbb{Z}^d(d\leq2)$ and transient in $\mathbb{Z}^d(d\geq3)$. In this context Holroyd and Propp(\cite{HOPR}) posed a question that for a Rotor walk on $\mathbb{Z}^d$ with $d\geq3$, does there exist a rotor configuration which is recurrent.

\section{Recurrent Rotor Configuration in $\mathbb{Z}^d$}
\paragraph{}
Here in this section we would like to suggest an initial Rotor configuration so that it exhibits the recurrence property. We will propose a rotor configuration in such a way that whenever the particle visits any site for the first time will be sent towards the origin. Mathematically it can be defined as follows.
\paragraph{}
Let $(x_1,x_2,...,x_d)$ be any point in $\mathbb{Z}^d$. Let $e_0,e_2,...,e_{d-1}$ be the unit vectors along the positive $X_1,X_2,...,X_d$ axes respectively and $e_{d},e_{d+1},...,e_{2d-1}$ be the unit vectors along the negative $X_1,X_2,...,X_d$ axes respectively. Now define the initial rotor configuration as follows:

\begin{equation*}
    \l(x_1,x_2,...,x_d) = \left\{
    \begin{array}{cc}
    e_{2d-1} &\quad \mbox \quad \mbox{if $|x_1|>max\{|x_2|,|x_3|,...,|x_d|\}$ and $x_1<0$}\\
    e_1 &\quad \mbox \quad \mbox{if $|x_2|>max\{|x_1|,|x_3|,...,|x_d|
    \}$ and $x_2<0$} \\
    \vdots &\quad \mbox \quad \mbox\quad \mbox \quad \mbox\quad \mbox \quad \mbox\vdots\\
    e_{d-1} &\quad \mbox \quad \mbox{if $|x_d|>max\{|x_1|,|x_2|,...,|x_{d-1}|\}$ and $x_d<0$}\\
    e_d &\quad \mbox \quad \mbox{if $|x_1|>max\{|x_2|,|x_3|,...,|x_d|
    \}$ and $x_1>0$} \\
    e_{d+1} &\quad \mbox \quad \mbox{if $|x_2|>max\{|x_1|,|x_3|,...,|x_d|\}$ and $x_2>0$}\\
    \vdots &
    \quad \mbox \quad \mbox\quad \mbox \quad \mbox\quad \mbox \quad \mbox\vdots\\
    e_{2d-2} &\quad \mbox \quad \mbox{if $|x_d|>max\{|x_1|,|x_2|,...,|x_{d-1}|\}$ and $x_d>0$}\\
    e_1 &\quad \mbox{in all other possible cases}
    \end{array} \right.
\end{equation*}

Also the sequence of the directions in which they rotate is ($e_0,e_1,...,e_{d-1},e_{d},e_{d+1},...,e_{2d-1}$).\\
The above will give us the initial rotor configuration in $\mathbb{Z}^d$.\\

 We claim that the above described rotor configuration is recurrent.
\subsection*{The case d=3}
\paragraph{}
For the case $d=3$, we have a strong numerical evidence to show that the above claim we made is true. Define the box to be $\mathbb{B}[0,n] = [-n,n]^d$. It is observed that the particle is visiting the origin 6n+1 times before it first leaves the box $\mathbb{B}[0,n]$. This has been checked for all $n\leq250$. It is also observed that for any point after some time, the direction assigned to it is same whenever the particle exits the box $\mathbb{B}[-n,n]$. The above mentioned facts clearly support our claim.

\section{What Next?}
\paragraph{}
If our claim gets proved, then one can attempt to classify the rotor configurations into transient and recurrent classes. Attempts can also be made to characterize these configuration by a simple rules. Studying these classes separately in detail may yield better results and can used as a replacement for classical probability models involving random walks. The transient configurations in $\mathbb{Z}^d(d\leq2)$ and recurrent configurations in $\mathbb{Z}^d(d\geq3)$ are the analogues of measure-zero sets in classical random walk. Hence the study about these in detail can be used to model the conditional events of these measure-zero sets.

\section{Acknowledgements}
\paragraph{}
I thank Dr.Arni S.R. Srinivasa Rao for exposing me to this emerging field, and also for guiding me in this project for the past few months. I also thank all of my friends who helped me in this project.

\end{document}